\documentclass[10pt]{amsart}

\usepackage{amssymb}

\newcommand{\F}{\mathcal{F}}

\newtheorem{theorem}{Theorem}[section]

\theoremstyle{definition}

\theoremstyle{remark}

\numberwithin{equation}{section}

\begin{document}

\title{Towards an automation of the circle method}


 \author[A. V. Sills]{Andrew V. Sills}
\address{Department of Mathematical Sciences, Georgia Southern University,
Statesboro, GA, 31407-8093, USA}
\email{ASills@GeorgiaSouthern.edu}



\date{\today}

\begin{abstract}
The derivation of the Hardy-Ramanujan-Rademacher formula for the
number of partitions of $n$ is reviewed.  Next, the steps for finding
analogous formulas for certain restricted classes of partitions or overpartiitons
is examined, bearing in mind how these calculations can be
automated in a CAS.  Finally, a number of new formulas of this type
which were conjectured with the aid of \emph{Mathematica}
are presented along with results of a test for their numerical accuracy.
\end{abstract}
\allowdisplaybreaks

\maketitle


\vskip 3cm

\section{Introduction}
A \emph{partition} of an integer $n$ is a representation of $n$ as a sum of positive integers,
where the order of the summands (called \emph{parts}) is considered irrelevant.
For example, there are seven partitions of the integer $5$, namely $5$, $4+1$, $3+2$,
$3+1+1$, $2+2+1$, $2+1+1+1$, and $1+1+1+1+1$.
Euler~\cite{le1748} was the first to systematically study partitions.  He showed that
\begin{equation}\label{egf}
\sum_{n=0}^\infty p(n) q^n = \prod_{m=1}^\infty \frac{1}{1-q^m} ,
\end{equation}
where $p(n)$ denotes the number of partitions of $n$ and we follow the convention
that $p(0)=1$.
  
   The series and infinite product in~\eqref{egf} converge absolutely when $|q|<1$.
Hardy and Ramanujan were the first to study $p(n)$ analytically and showed that~\cite[p. 79, Eq. (1.41)]{hr1918}
\begin{equation} \label{hrasymp}
p(n) \sim \frac{ \exp (\pi\sqrt{2/3}) }{4n\sqrt{3}} \qquad \mbox{as $n\to\infty$}.
\end{equation}
Noting that the value of $p(200)$ estimated by~\eqref{hrasymp} was surprisingly
close to the true value of $p(200)$ as computed by P.A. MacMahon, Hardy and Ramanujan
were encouraged to push their analysis of $p(n)$ further.
  
  Ultimately, they produced the formula~\cite[p. 85, Eq. (1.75)]{hr1918}
\begin{multline} \label{hr} p(n) = \frac{1}{2\pi\sqrt{2}} \sum_{k=1}^{\alpha\sqrt{n}} \sqrt{k} 
\underset{(h,k)=1}{\sum_{0\leqq h < k}} \omega(h,k) e^{-2\pi i h n/k}
\frac{d}{dn} \left( \frac{\exp\left( \frac{\pi}{k} \sqrt{\frac 23 (n-\frac{1}{24} ) }  \right)}{\sqrt{n-\frac{1}{24}}} \right)  \\ + O(n^{-1/4}),
\end{multline}
with $\alpha$ an arbitrary constant and $\omega(h,k)$ a certain complex $24k$th root of 
unity\footnote{Apostol~\cite{tma1948} showed that it is also a $12k$th root of unity.}
that arises frequently in the theory of modular forms and is defined below
in~\eqref{OmegaDef}.

  While~\eqref{hr} is an asymptotic formula, it is incredibly accurate.  For the case
$n=200$, summing $k$ from $1$ to $8$ results in a value only $0.004$ higher than
the true value of $3\ 972\ 999\ 029\ 338$.

Later, it was shown by D. H. Lehmer~\cite{dhl1937} that if the sum on $k$ in~\eqref{hr}
is extended to $\infty$, the resulting series diverges.

In~\cite{har1937}, Rademacher made a slight change in Hardy and
Ramanujan's analysis which led him to finding a convergent series
representation for $p(n)$, very similar in form to that of~\eqref{hr}.
This result is presented below as Theorem~\ref{RadPofN}.
In a later paper,
Rademacher altered the path of integration and as a result was able to give a simpler
proof for the correctness of his series~\cite{har1943}.  This latter technique is also
described in books by Rademacher~\cite[Ch. 14]{har1973} and Apostol~\cite[Ch. 5]{tma1990},
while the former may be found in the text of Andrews~\cite[Ch. 5]{gea1976}.

The technique of deriving the formula for $p(n)$ via integration of 
a certain function (see \eqref{Cauchy} below) which has singularities at every point of the unit circle
in the complex plane has come to be known as the ``circle method."
The circle method has proven to be applicable to many problems
and as such is one of the most
important and useful tools in analytic number theory.  
There are far too many papers which have used the circle method to
even begin to mention them here, but a subset of the literature which
employs the circle method to find formulas for certain restricted classes
of partitions includes Grosswald~\cite{eg1958, eg1960}, Haberzetle~\cite{mh1941},
Hagis~\cite{ph1962, ph1963, ph1964a, ph1964b, ph1965a, ph1965b, ph1965c, ph1966, ph1971}, Hua~\cite{lkh1942}, Iseki~\cite{si1959, si1960, si1961}, Lehner~\cite{jl1941}, Livingood~\cite{jl1945}, 
Niven~\cite{in1940}, and Subramanyasastri~\cite{vvs1972}.
Recently, Bringmann and Ono~\cite{bo2009} have made great strides
in the subject by finding exact formulas for the coefficients of harmonic
Maas forms of nonpositive weight.

   A main theme of this paper is that while the application of the circle
method to find $p(n)$ or a given restricted partition formula may be
complicated, it is essentially a \emph{calculation}. 
As such, many of the steps involved are ripe for
automation.  
Furthermore, a good number of the steps involve showing that a given integral approaches zero.  As long as we can 
reliably predict when this will be the case, we can produce 
reasonable conjectures for formulas without worrying
about the estimates that are required when a rigorous proof is desired.  
 
 We shall outline a derivation of $p(n)$, and then consider how the
circle method applies to
restricted partition and overpartition formulas, bearing in mind how
to automate these calculations.

 Finally, we shall present some new restricted formulas conjectured with
the aid of \emph{Mathematica}.
\pagebreak

 \section[Hardy-Ramanujan-Rademacher Formula for $p(n)$]{An Overview of the Derivation of the\\ Hardy-Ramanujan-Rademacher Formula for $p(n)$}
\subsection{Preliminaries}
\subsubsection{The Dedekind $\eta$-function}
Let $\mathcal H:= \{ \tau\in\mathbb{C} \ | \ \Im\tau>0 \},$ the upper half of the complex plane. 

The Dedekind eta function is defined by
\begin{equation} \label{EtaDef}
\eta(\tau) := e^{\pi i \tau/ 12} \prod_{m=1}^\infty (1- e^{2\pi i m \tau})
\end{equation}
where $\tau\in\mathcal{H}$.

 For $a$, $b$, $c$, $d \in \mathbb{Z}$ with $ad-bc=0$, and $c>0$,
 $\eta(\tau)$ satisfies the
 functional equation
 \begin{equation} \label{EtaFunctionalEq}
 \eta\left( \frac{a\tau+b}{c\tau +d} \right)=
  \omega(-d, c)  \exp\Bigg( \pi i \left( \frac{a+d}{12c} \right) \Bigg)
  \sqrt{-i(c\tau +d)}\  \eta(\tau),
 \end{equation}
 where 
 \begin{equation}\label{OmegaDef}
 \omega(h,k) = \left\{  
   \begin{array}{ll}
       \left( \frac{-k}{h} \right) \exp \left(  -\pi i \left \{ \frac 14 (2-hk-h) + \frac{1}{12}(k - \frac 1k)
       (2h-H+h^2 H) \right\} \right), &\mbox{if $2\nmid h$}\\
       \left( \frac{-h}{k} \right) \exp \left( -\pi i \left\{ \frac 14 (k-1) + \frac{1}{12}(k- \frac 1k)(2h-H_h^2 H) \right\} \right),  &\mbox{if $2\nmid k$} \\
   \end{array}
 \right.
\end{equation}
$( \frac ab )$ is the Legendre-Jacobi symbol, and $H$ is any solution of the
congruence
\[ hH \equiv -1 \pmod{k}. \]


\subsubsection{Farey fractions}

The sequence $\F_N$ of \emph{proper Farey fractions of order $N$} is the set of all $\frac hk$ with
$(h,k)=1$ and $0\leqq \frac hk <1$, arranged in increasing order.

Thus, we have
\begin{equation*}
\F_1 = \left\{ \frac 01 \right\}, \quad
\F_2 = \left\{ \frac 01, \frac 12 \right\}, \quad
\F_3 = \left\{ \frac 01, \frac 13, \frac 12, \frac 23 \right\}, \quad
\F_4 = \left\{ \frac 01, \frac 14, \frac 13, \frac 12, \frac 23, \frac 34\right\},
\end{equation*}
etc.

For a given $N$, let $h_p$, $h_s$, $k_p$, and $k_s$ be such that $\frac{h_p}{k_p}$ is
the immediate predecessor of $\frac hk $ and $\frac{h_s}{k_s}$ is the immediate
successor of $\frac hk$ in $\F_N$.  It will be convenient to view each $F_N$ cyclically, i.e.
to view
 $\frac 01$ as the immediate successor of
$\frac {N-1}{N}$.


\subsubsection{Ford circles and the Rademacher path}
Let $h$ and $k$ be integers with $(h,k) = 1$ and $0\leqq h < k$.  
The \emph{Ford circle}~\cite{lrf1938}
$C(h,k)$ is the circle in $\mathbb C$ of radius $\frac{1}{2k^2}$ centered at the point
$$ \frac{h}{k} + \frac{1}{2k^2} i.$$  

The \emph{upper arc $\gamma(h,k)$ of the Ford circle $C(h,k)$} 
is the arc of the circle
\[ \left| \tau - \left(\frac hk + \frac{1}{2k^2} i \right) \right| = \frac{1}{2k} \]
from the initial point 
\begin{equation} \label{AlphaI}
 \alpha_I(h,k):=  \frac hk - \frac{k_p}{k(k^2+k_p^2)} + \frac{1}{k^2+k_p^2} i
\end{equation}
to the terminal point 
\begin{equation}\label{AlphaT}
 \alpha_T(h,k):= \frac hk + \frac{k_s}{k(k^2+k_s^2)} + \frac{1}{k^2+k_s^2} i, 
 \end{equation}
traversed clockwise.

Note that we have
$\alpha_I(0,1)  = \alpha_T(N-1,N) . $

Every Ford circle is in the upper half plane.
For $\frac{h_1}{k_1}, \frac{h_2}{k_2} \in \F_N$,  $C(h_1, k_1)$ and $C(h_2, k_2)$ are
either tangent or do not intersect.

The \emph{Rademacher path} $P(N)$ of order $N$ is the path in the upper half of the
$\tau$-plane from $i$ to $i+1$ consisting of 
\begin{equation} \label{RadPath}  \bigcup_{\frac hk \in \F_N}  \gamma(h,k) \end{equation}
traversed left to right and clockwise.  In particular, we consider the left half of the Ford
circle $C(0,1)$ and the corresponding upper arc $\gamma(0,1)$ to be translated
to the right by 1 unit.  This is legal given then periodicity of the function which is to 
be integrated over $P(N)$.

\subsection{Euler and Cauchy get us off the ground}
Recall Euler's generating function for $p(n)$,
\begin{equation}
f(q):=  \sum_{n=0}^\infty p(n) q^n = \prod_{m=1}^\infty \frac{1}{1-q^m}.
\end{equation}

 Let us now fix $n$.
The function $f(q) / q^{n+1} $ has a pole of order $n+1$ at $q=0$, and
an essential singularity at every point of the unit circle $|q|=1$. 
The Laurent series of  $f(q) / q^{n+1} $ about $q=0$ is therefore
\[ \sum_{j=0}^\infty p(j) q^{j-n-1} = \sum_{j=-n-1}^\infty p(j+n+1) q^{j} , \]
for $0<|q|<1$,
and so the residue of $f(q) / q^{n+1}$ at $q=0$ is $p(n)$.

Thus, Cauchy's integral formula implies that
\begin{equation} \label{Cauchy}
p(n) = \frac{1}{2\pi i}\int_{\mathcal{C}} \frac{f(q)}{q^{n+1}} \ dq, 
\end{equation}
where $\mathcal{C}$ is any positively oriented, simple closed contour enclosing the origin and
inside the unit circle. 

\subsection{The choice of $\mathcal{C}$}

Since \begin{equation*}
\frac{f(q)}{q^{n+1}}  = \frac{1}{q^{n+1}}\prod_{k=1}^\infty \frac{1}{1-q^k}
= \frac{1}{q^{n+1}}\prod_{k=1}^\infty \prod_{j=0}^{k-1} \frac{1}{1- e^{2\pi i j/k}q}
\end{equation*}
we see that although every point of along $|q|=1$ is an essential singularity of $f(q)/ q^{n+1}$,
in some sense $q=1$ is the ``heaviest" singularity, $q=-1$ is ``half as heavy," 
$q= e^{2\pi i/3}$ and $e^{4 \pi i/3}$ are each ``one third as heavy," etc.

The integral~\eqref{Cauchy} is evaluated by approximating the integrand for each 
$h, k$ by an elementary function which is very nearly equal to $f(q) / q^{n+1}$ near
the singularity $e^{2\pi i h/k}$.  The contor $\mathcal C$ is chosen in such a way that
the error introduced by this approximation is carefully kept under control.

  We introduce the change of variable 
$ q = \exp(2 \pi i \tau) $
so that the unit disk $|q|\leqq 1$ in the $q$-plane maps to the infinitely tall, unit wide strip in the
$\tau$ plane where $0\leqq \Re\tau\leqq 1$ and $\Im\tau\geqq 0$.  The contour $\mathcal{C}$
is then taken to be the preimage of the Rademacher path $P(N)$  (see~\eqref{RadPath})
under the map 
$q \mapsto \exp(2 \pi i \tau) $.  Better yet, let us replace $q$ with 
$\exp(2 \pi i \tau)$ in~\eqref{Cauchy}
to express the integration in the $\tau$-plane:
\begin{align*}\label{TauIntegral}
p(n) &= \int_{P(N)} f(e^{2\pi i \tau}) e^{-2\pi i n\tau} d\tau\\
        &=\sum_{\frac hk \in \F_N} \int_{\gamma(h,k)} f(e^{2\pi i \tau}) e^{-2\pi i n \tau} d\tau\\
         &=\sum_{k=1}^{N} \underset{(h,k)=1}{ \sum_{0 \leqq h < k}} \int_{\gamma(h,k)} 
         f(e^{2\pi i \tau}) e^{-2\pi i n \tau} d\tau
\end{align*}

\subsection{Another change of variable}
Next, we change variables again, taking
\begin{equation}\label{TauToZ} \tau = \frac {iz + h}{k}\end{equation} so that
$ z = -ik \left(\tau - \frac hk\right),$
for each $\tau \in C(h,k)$.
Thus $C(h,k)$ (in the $\tau$-plane) maps to the clockwise-oriented 
circle $K^{(-)}_k$ (in the $z$-plane) centered
at ${1}/{2k}$ with radius ${1}/{2k}$.

So we now have 
\begin{align}
p(n) &= \sum_{k=1}^{N} \underset{(h,k)=1}
{ \sum_{0 \leqq h < k}} \int_{z_I(h,k)}^{z_T(h,k)} 
f( e^{2\pi i h/k - 2\pi z/k} ) e^{-2\pi i n (iz + h)/k} \frac{i}{k}\  dz \\
&= \sum_{k=1}^{N} \underset{(h,k)=1}
{ \sum_{0 \leqq h <  k}} \frac ik \ e^{-2\pi i n h/k} 
\int_{z_I(h,k)}^{z_T(h,k)} e^{2n\pi z/k} f( e^{2\pi i h/k - 2\pi z/k} ) \ dz, \label{IntegralZ}
\end{align}
where $z_I(h,k)$ (resp. $z_T(h,k)$) is the image of $\alpha_I(h,k)$ (see~\eqref{AlphaI})
(resp. $\alpha_T(h,k)$ [see~\eqref{AlphaT}]) under the transformation~\eqref{TauToZ}.

So the transformation~\eqref{TauToZ} maps the upper arc $\gamma(h,k)$ of $C(h,k)$
in the $\tau$-plane to the arc on $K^{(-)}_k$ which initiates at
\begin{equation} \label{ZI}
 z_I(h,k) = \frac{k}{k^2 + k_p^2} + \frac{ k_p}{k^2+ k_p^2} i
\end{equation} and terminates at
\begin{equation} \label{ZT}
  z_T(h,k) = \frac{k}{k^2+k_s^2} - \frac{ k_s}{k^2+k_s^2} i.
\end{equation}

\subsection{Exploiting a modular transformation}
It is incredibly fortunate that 
\[ f(q) = f(e^{2\pi i \tau}) = \frac{e^{\pi i \tau/12} }{\eta(\tau)}, \]
so that we may take advantage of the modular functional equation~\eqref{EtaFunctionalEq}
satisfied by $\eta(\tau)$ in our effort to evaluate~\eqref{IntegralZ}.
Equation~\eqref{EtaFunctionalEq} rewritten in terms of $f(q)$ is
\begin{equation} \label{FFunctionalEq}
 f(e^{2\pi i(iz+h)/k}) = \omega(h,k) e^{\pi(z^{-1}-z)/12k} \sqrt{z} f( e^{2\pi i( iz^{-1}+H)/k} ),
\end{equation}
where $(h,k)=1$ and $H$ is a solution to the congruence
$hH \equiv -1 \pmod{k} . $

 Note that when $|z|$ is close to $0$, the left hand side of~\eqref{FFunctionalEq} is close
to $f( e^{2\pi i h/k})$, i.e. for $|z|$ small,~\eqref{FFunctionalEq} gives a good approximation
for $f$ evaluated at the ``heavy" singularity $e^{2 \pi i h/k}$.  Next, observe that the
final factor on the right hand side of~\eqref{FFunctionalEq},
\[ f(e^{2\pi i(iz^{-1}+H)/k}) = f\Bigg(\exp \left( \frac{2\pi i H}{k} - \frac{2\pi}{zk} \right) \Bigg), \]
is close to $ f(0) = 1$ when $|z|$ is small, so that
\[ f(e^{2\pi i(iz^{-1}+H)/k}) -1 \]
is close to $0$ when $|z|$ is small.

Applying this information to~\eqref{IntegralZ}, we find that
\begin{align}
p(n) 
&= \sum_{k=1}^{N} \underset{(h,k)=1}
{ \sum_{0 \leqq h < k}} \frac ik \ e^{-2\pi i n h/k} \omega(h,k)\nonumber \\
& \qquad \times
\int_{z_I(h,k)}^{z_T(h,k)} e^{2n\pi z/k }
e^{ \pi(z^{-1}-z)/12k} \sqrt{z} f( e^{2\pi i( iz^{-1}+H)/k} ) \ dz \nonumber \\
&= \sum_{k=1}^{N} \frac ik \underset{(h,k)=1}
{ \sum_{0 \leqq h < k}}  \ e^{-2\pi i n h/k} \omega(h,k) \nonumber \\ &\qquad \times
\int_{z_I(h,k)}^{z_T(h,k)} 
e^{ \pi (24 nz + z^{-1} - z)/12k} \sqrt{z} \left\{ 1+ \left[ f( e^{2\pi i( iz^{-1}+H)/k} ) -1\right]
\right\} \ dz 
\label{TransIntegralZ}\\
&= \sum_{k=1}^{N} \frac ik \underset{(h,k)=1}
{ \sum_{0 \leqq h < k}}  \ e^{-2\pi i n h/k} \omega(h,k) \left(I_{h,k} + I^*_{h,k}\right), \nonumber
\end{align}
where
\begin{equation} \label{I}
I_{h,k}:= \int_{z_I(h,k)}^{z_T(h,k)} 
e^{ \pi (24 nz + z^{-1} - z)/12k} \sqrt{z} \ dz 
\end{equation}
and
\begin{equation}\label{Istar}
I^*_{h,k}:=
\int_{z_I(h,k)}^{z_T(h,k)} 
e^{ \pi (24 nz + z^{-1} - z)/12k} \sqrt{z} \left[ f( e^{2\pi i( iz^{-1}+H)/k} ) -1\right] \ dz .
\end{equation}

\subsection{Estimating $I^*_{h,k}$}
The next goal is to show that $I^*_{h,k}$ is small when $N$ is large.
Note that we can change the path of integration of~\eqref{Istar} from an arc of the circle 
that is the image of the Ford circle under the transformation~\eqref{TauToZ} connecting
$z_I(h,k)$ and $z_T(h,k)$ to the line segment connecting $z_I(h,k)$ and $z_T(h,k)$
without altering the value of the integral.  On the segment connecting $z_I(h,k)$ and
$z_T(h,k)$, we have
\begin{equation} \label{AbsZMax}
 |z| \leqq \max \left\{ |z_I(h,k)|, |z_T(h,k)| \right\} = 
 \max \left\{  \sqrt{ \frac{1}{k^2+k_p^2} } ,\sqrt{ \frac{1}{ k^2+k_s^2 } }  \right\}
\leqq 
  \frac{\sqrt{2}}{N}  
 \end{equation}
 
Obviously, the length of the segment connecting $z_I(h,k)$ and $z_T(h,k)$ can be
easily calculated for any particular $h$ and $k$.  However, we wish to have an 
upper bound for the length that holds for a given $N$. 

The length of the segment is
$\leqq |z_I(h,k)| + |z_T(h,k)| \leqq \frac{2\sqrt{2}}{N}$.
 
 Bearing in mind that on the segment, 
 $\Re z < \frac 1k \mbox{ and } \Re \left( \frac 1z \right) = k, $
it can be shown that the integrand in~\eqref{Istar} is less than
$c|z|^{1/2}$, where
\[ c = e^{2n\pi/k^2} \sum_{m=1}^\infty p(24m-1) t^{24m-1} ,\]
with $t=e^{-\frac{\pi}{12}}$,
by mimicking the argument in~\cite[p. 107]{tma1990}.
 

Since $z$ is on the segment connecting $z_I(h,k)$ to $z_T(h,k)$, $|z|$ is bounded above
by $\sqrt{2}/N$, so the integrand
is bounded above by $c \sqrt{ \sqrt{2}/N }  = c 2^{1/4} N^{-1/2}$.
  Thus
 \[ | I^*_{h,k} | <   \frac{c 2^{1/4}}{N^{1/2}} \frac{2\sqrt{2}}{N}  =  C N^{-3/2}, \]
 where $C =  2^{7/4} c$.
 
 Finally, it can be shown that 
 $\left|\sum_{k=1}^N\frac ik\sum_{0\le h<k,\,(h,k)=1} e^{-2\pi inh/k}\omega(h,k)I^*_{h,k}\right|=O(N^{-1/2})$
 
 \subsection{Estimations associated with $I_{h,k}$}
 The work of the preceding section allows us to rewrite~\eqref{TransIntegralZ} as
 \begin{equation}
 p(n)= \sum_{k=1}^{N} \frac ik \underset{(h,k)=1}
{ \sum_{0 \leqq h < k}}  \ e^{-2\pi i n h/k} \omega(h,k) I_{h,k} + O(N^{-1/2}),
\end{equation} where, as before,
\begin{equation} \label{Irestate}
I_{h,k}:= \int_{z_I(h,k)}^{z_T(h,k)} 
e^{ \pi (24 nz + z^{-1} - z)/12k} \sqrt{z} \ dz .
\end{equation}
We proceed by re\"expressing $I_{h,k}$ as
\begin{equation}
 I_{h,k} = \int_{K^{(-)}_k} - \int_0^{z_I(h,k)} - \int_{z_T(h,k)}^0,
\end{equation}
where the integrands of all three integrals 
are the same as that of the right hand side of~\eqref{Irestate}.

  The length of the arc connecting $0$ and $z_I(h,k)$ is less than
  \[ \frac{\pi}{2} | z_I(h,k) | < \frac{\pi}{2} \frac{ \sqrt{2}}{N}. \]
 On the arc, $|z| < \sqrt{2}/N$.
 
   We had previously seen that the absolute value of the integrand is $< c |z|^{1/2}$,
 so
 \[ \int_0^{z_I(h,k)} e^{ \pi (24 nz + z^{-1} - z)/12k} \sqrt{z} \ dz 
 < c \sqrt{ \frac{\sqrt{2}}{N} }\frac{\pi}{\sqrt{2} N} = C N^{-3/2} .\]
 An analogous estimate applies to $\int_{z_T(h,)}^0 $.
 
 \subsection{The formula for $p(n)$}
 We may now write
 \begin{equation}
p(n) = \sum_{k=1}^{N} \underset{(h,k)=1} 
{ \sum_{0 \leqq h < k}} \frac ik \ e^{-2\pi i n h/k} 
\int_{K^{(-)}_k} e^{2n\pi z/k } \omega(h,k)
e^{ \pi(z^{-1}-z)/12k} \sqrt{z} \ dz + O(N^{-1/2}).
 \end{equation}
Let $N\to\infty$ to obtain
 \begin{equation}
p(n) = i\sum_{k=1}^{\infty} \underset{(h,k)=1}
{ \sum_{0 \leqq h < k}} \frac{ e^{-2\pi i n h/k}}{k} \omega(h,k) 
\int_{K^{(-)}_k} \sqrt{z}  
\exp \left\{ \frac{\pi}{12zk} + \frac{2\pi z}{k} \left( n - \frac{1}{24} \right) \right\} \ dz .
 \end{equation}

Next, apply the transformation $z=1/w$ so that $dz = -1/w^2\ dw$:
 \begin{equation} \label{WIntegral}
p(n) =\frac{1}{i}\sum_{k=1}^{\infty} \frac 1k\underset{(h,k)=1}
{ \sum_{0 \leqq h < k}}  \ e^{-2\pi i n h/k } \omega_{h,k} 
\int_{1-\infty i}^{1+\infty i} w^{-5/2}  
\exp \left\{ \frac{\pi w}{12k} + \frac{2\pi}{wk} \left( n - \frac{1}{24} \right) \right\} \ dw .
 \end{equation}
 The integral in~\eqref{WIntegral} can be evaluated in terms of Bessel functions.
To make this evaluation easier to see, we set 
$ w= \frac{12k}{\pi} t $, so that $dw = \frac{12k}{\pi} dt$, to obtain
  \begin{multline} \label{TIntegral}
p(n) =2\pi \left(\frac{\pi}{12}\right)^{3/2}\sum_{k=1}^{\infty} \underset{(h,k)=1}
{ \sum_{0 \leqq h < k}} k^{-5/2} \ e^{-2\pi i n h/k } \omega(h,k) \\ \times
\frac{1}{2\pi i}\int_{\pi/12k -\infty i}^{\pi/12k +\infty i} t^{-5/2}  
\exp \left\{ t + \frac{\pi^2}{6 k^2 t} \left( n - \frac{1}{24} \right) \right\} \ dt .
 \end{multline}
 
Now recall the Bessel function of the first kind of purely imaginary
argument is given by~\cite[p. 181, Eq. (1)]{gnw1944}
 \[ I_\nu(z)= \frac{(z/2)^\nu}{2\pi i} \int_{-\infty}^{(0+)} t^{-\nu-1} \exp \left( t + \frac{z^2}{4t} \right) \ dt.\]
Taking into account the remark preceding Eq. (8) on p. 177 of~\cite{gnw1944}, we may, since
$\pi/12k > 0$, alter
the path of integration to obtain
 \begin{equation} \label{BesselI}
  I_\nu(z)= \frac{(z/2)^\nu}{2\pi i} \int_{\pi/12k-i\infty}^{\pi/12k+i\infty} t^{-\nu-1} \exp \left( t + \frac{z^2}{4t} \right) \ dt.
  \end{equation}
 Setting $\nu=3/2$ and $z=\frac{\pi}{k}\sqrt{ \frac 23 (n-\frac{1}{24}) } $ in~\eqref{BesselI}
 and applying the result to~\eqref{TIntegral}, we find
 \begin{equation}
 p(n) = \frac{2\pi}{(24)^{3/2}} \left( n- \frac{1}{24} \right)^{-3/4} \sum_{k=1}^\infty \frac 1k
 \underset{(h,k)=1}
{ \sum_{0 \leqq h < k}}  e^{ - 2\pi i nh/k} \omega_{h,k} I_{3/2} \left( \frac{\pi}{k}
 \sqrt{\frac 23 \left(n-\frac{1}{24}\right)}
\right)
\end{equation}
Bessel functions of half-odd order can be written in terms of elementary functions.
In particular,
  \[ I_{3/2}(z) = \sqrt{\frac{2z}{\pi}} \frac{d}{dz} \left( \frac{ \sinh z}{z} \right), \]
so the final form of the formula for $p(n)$ is
\begin{theorem}[Rademacher] \label{RadPofN}
\begin{equation*}
p(n) = \frac{1}{\pi\sqrt{2}} \sum_{k=1}^\infty \sqrt{k}
 \underset{(h,k)=1}
{ \sum_{0 \leqq h < k}}  e^{ - 2\pi i nh/k} \omega(h,k)
\frac{d}{dn} \left( \frac{ \sinh \left( \frac{\pi}{k}  \sqrt{\frac23 \left( n-\frac{1}{24}\right)}\right)}{\sqrt{n-\frac{1}{24}}} \right).
\end{equation*} 
 \end{theorem}

\section{Restricted Partition Functions} 
\subsection{Partition Identities}
Euler~\cite{le1748} observed that the algebraic identity
\begin{equation} \label{OddDistinctGF}
\prod_{j=1}^\infty (1+q^j) = \prod_{j=1}^\infty \frac{1}{1-q^{2j-1}}
\end{equation}
implies the following theorem about integer partitions:
\begin{theorem}[Euler]
The number of partitions of $n$ into distinct parts equals the number of partitions of $n$
into odd parts.
\end{theorem}

While such a result tells us that there are the same number of partitions of
$n$ into
distinct parts as there are partitions of $n$ using only odd parts,
we do not know 
\emph{how many} such partitions of $n$ there are.  The circle method has
been applied by P. Hagis~\cite{ph1963} and~L. K. Hua~\cite{lkh1942} to address this question.

\begin{theorem}[Hagis] \label{Hagis1} Let $\delta(n)$ denote the
number of partitions of $n$ into distinct parts.  Then
\begin{equation} \label{HagisDist}
 \delta(n) = \frac{\pi}{\sqrt{24n+1}}
 \underset{2\nmid k}{ \sum_{k\geqq 1}} \frac 1k
  \underset{(h,k)=1}{\sum_{0\leqq h <k}}
    e^{-2\pi n h /k}\frac{\omega(h,k)}{\omega(2h,k)} I_1\left( \frac{\pi\sqrt{24n+1}}{6\sqrt{2} k} \right).
    \end{equation}
\end{theorem}

J.W.L. Glaisher~\cite{jwlg1883} generalized Euler's result to
\begin{theorem}[Glaisher]
The number of partitions of $n$ where no part appears more than $j-1$ times equals
the number of ``$j$-regular partitions of $n$", i.e.
partitions of $n$ where no part is a multiple of $j$.
\end{theorem} 
Clearly,
Euler's theorem is the $j=2$ case of Glaisher's theorem.  Glaisher's theorem follows
immediately from the identity
\begin{equation} \label{GlaisherGF}
\prod_{k=1}^\infty  (1+q^k + q^{2k} + \cdots + q^{(j-1)k} )
= \underset{k\not\equiv 0\hskip-3mm\pmod{j}}{\prod_{k\geqq 1}} \frac{1}{1-q^k}.
\end{equation}

\begin{theorem}[Hagis~\cite{ph1971}] \label{Hagis2}
Let $\delta_j (n)$ denote the number of $j$-regular partitions of $n$.  Then
\begin{multline} \delta_j(n) = \frac{2 \pi}{j\sqrt{24n+j-1}}
  \underset{d\mid j}{\sum_{0<d<\sqrt{j}}} \sqrt{d(j-d^2) }
 \underset{(k,j)=d}{\sum_{k\geqq 1}}
  \frac 1 k \\ \times
 \underset{(h,k)=1}{\sum_{0\leqq h <k}}
    e^{-2\pi n h /k}\frac{\omega(h,k)}{ \omega( \frac{jh}{d} , \frac kd)  }
     I_1\left(  \frac{\pi}{6k} \sqrt{ \frac{(24n+j-1)(j-d^2)}{j} } \right).
   \end{multline}
\end{theorem}


  
Another well known partition identity of this type is
\begin{theorem}[I. Schur~\cite{is1926}]
The number of partitions of $n$ into distinct parts which differ by at least three and where no
consecutive multiples of three appear equals the number of partitions of $n$ into parts
congruent to $\pm 1\pmod{6}$.
\end{theorem}

\begin{theorem}[I. Niven~\cite{in1940}]  
Let $S(n)$ denote the number of partitions of
$n$ into parts congruent to $\pm 1\pmod{6}$.  Then
   \begin{multline}
  S(n) = \frac{\pi}{\sqrt{36n-3}} \sum_{d\mid 6} \sqrt{(d-2)(d-3)}
    \underset{(k,6)=d}{\sum_{k\geqq 1}} \frac 1k 
   \\ \times \underset{(h,k)=1}{\sum_{0\leqq h <k}} 
    e^{-2\pi n h /k}
    \frac{\omega(h,k) \omega( \frac{6h}{(k,6)}, \frac{k}{(k,6)} )}
    {\omega(\frac{2h}{(k,2)},\frac{k}{(k,2)})   \omega(\frac{3h}{(k,3)},  \frac{k}{(k,3)})} 
    I_1\left( \frac{\pi\sqrt{d\left(12n-1\right)}}{3 \sqrt{6} k} \right).
    \label{Niven} 
   \end{multline}
\end{theorem}

Recently, the author found~\cite{avs2009}
 \begin{equation}
 \label{overptnformula}
  \bar{p}(n) = \frac{1}{2\pi} \underset{2\nmid k}{\sum_{k\geqq 1}} \sqrt{k}
 \underset{(h,k)=1}{\sum_{0\leqq h < k}} \frac{\omega(h,k)^2 }{ \omega(2h,k)  }
     e^{-2\pi i nh /k } \frac{d}{dn} \left(  \frac{\sinh \left( \frac{\pi \sqrt{n}}{k} \right) }{\sqrt{n} } \right)
 \end{equation} and
   \begin{multline} pod(n) = \frac{2}{\pi\sqrt{6}}
   \sum_{d\mid 4} \sqrt{(d-2)(5d-17)}
    \underset{(k,4)=d }{\sum_{k\geqq 1}}
 \sqrt{k}  \\ \times
 \underset{(h,k)=1}{\sum_{0\leqq h < k}} \frac{\omega(h,k) \ \omega\left(  \frac{4h}{d}, \frac{k}{d} \right) }{ \omega\left( \frac{2h}{(k,2)} ,\frac{k}{(k,2)} \right)  }
     e^{-2\pi i nh /k } \frac{d}{dn} \left(  \frac{\sinh \left( \frac{\pi \sqrt{d (8n-1)}}{4k} \right) }{\sqrt{8n-1} } \right),  \label{PodOfN}
 \end{multline}
where $pod(n)$ denotes the number of partitions of $n$ where no
odd part is repeated, and $\bar{p}(n)$ denotes the number of
overpartitions of $n$.  An~\emph{overpartition} of $n$ is a 
finite weakly decreasing sequence of positive integers
where the last occurrence of a given part may or may not be overlined.  Thus the eight overpartitions of $3$ are $(3)$, $(\bar{3})$, $(2,1)$, 
$(\bar{2},1)$, $(2,\bar{1})$, $(\bar{2},\bar{1})$, 
$(1,1,1)$, $(1,1,\bar{1})$.   Overpartitions were introduced by S. Corteel and J. Lovejoy in~\cite{cl2004} and have been studied extensively by them and others.

\subsection{Distinct Parts}

We have
\[ \sum_{n=0}^\infty \delta (n) q^n = \underset{2\nmid  j}{\prod_{j\geqq 1} }\frac{1}{1-q^j} = 
\frac{f(q)}{f(q^2)} =: F(q),\]
where, as before $f(q):= \sum_{n\geqq 0} p(n) q^n = \prod_{j\geqq 1} (1-q^j)^{-1}$.

 Proceeding as in the case of $p(n)$, we note
\begin{align}  
   \delta(n) 
   &= \frac{1}{2\pi i} \int_{\mathcal C} \frac{F(q)}{q^{n+1}} \ dq  \nonumber\\
   &= \frac{1}{2\pi i} \int_{\mathcal C} \frac{f(q)}{f(q^2) q^{n+1}} \ dq \nonumber\\
   &= \int_{P(N)}\frac{ f(e^{2\pi i \tau}) }{ f (e^{4\pi i \tau})}e^{-2\pi i n\tau} d\tau\nonumber\\
        &=\sum_{\frac hk \in \F_N} \int_{\gamma(h,k)} \frac{ f(e^{2\pi i \tau}) }{ f (e^{4\pi i \tau})} e^{-2\pi i n \tau} d\tau\nonumber\\
         &=\sum_{k=1}^{N} \underset{(h,k)=1}{ \sum_{0 \leqq h < k}} \int_{\gamma(h,k)} 
         \frac{ f(e^{2\pi i \tau}) }{ f (e^{4\pi i \tau})} e^{-2\pi i n \tau} d\tau\nonumber \\
      &= \sum_{k=1}^{N} \underset{(h,k)=1}
{ \sum_{0 \leqq h < k}} \int_{z_I(h,k)}^{z_T(h,k)} 
\frac{ f( e^{2\pi i h/k - 2\pi z/k} )}{f( e^{4\pi i h/k - 4\pi z/k} )} e^{-2\pi i n (iz + h)/k} \frac{i}{k}\  dz 
\nonumber\\
&= i\sum_{k=1}^{N} k^{-1} \underset{(h,k)=1}
{ \sum_{0 \leqq h <  k}}  \ e^{-2\pi i n h/k} 
\int_{z_I(h,k)}^{z_T(h,k)} e^{2n\pi z/k} 
\frac{ f( e^{2\pi i h/k - 2\pi z/k} )}{f( e^{4\pi i h/k - 4\pi z/k} )} \ dz, \label{OddIntegralZ}
\end{align}
where, as before, $q= e^{2\pi i \tau}$, $\tau = (i z + h)/k$, and $P(N)$, $\gamma(h,k)$,
$z_I(h,k)$, and $z_T(h,k)$ all have the same meaning as before.

At this point, we should like to transform 
$$F(q) = f(q)/f(q^2) = { f( e^{2\pi i h/k - 2\pi z/k} )}/{f( e^{4\pi i h/k - 4\pi z/k} )},$$
just as we had transformed  $f(q) = f( e^{2\pi i h/k - 2\pi z/k} )$
via~\eqref{FFunctionalEq} in the analogous analysis of $p(n)$.

It will be necessary to consider two cases.  
When $k$ is even, $k/2$ is an integer, so
we can obtain $f(q^2)$ from $f(q)$ by replacing $k$ by $k/2$ in $f( e^{2\pi i h/k - 2\pi z/k} )$.
On the other hand, when $k$ is odd, we instead replace $h$ by $2h$ and $z$ by $2z$
in $f( e^{2\pi i h/k - 2\pi z/k} )$.  Thus,
\begin{equation} 
F( e^{2\pi i h/k - 2\pi z/k} ) = 
  \left\{ 
  \begin{array}{ll} 
  \frac{\omega(h,k)}{\omega(h,k/2)}
    \exp\left( 
        \frac{\pi z}{12k} - \frac{\pi}{12kz}  \right)
     F\left(  \exp \left( \frac{2\pi i (H_1 + iz^{-1})  }{k} \right) \right), &\mbox{if $2 \mid k$,}\\
  { \frac{\omega(h,k)}{\omega(2h,k) \sqrt{2}}  \exp\left( \frac{\pi z}{12k}  + \frac{\pi}{24 kz}  \right) }/
    {F\left(  \exp \left( \frac{ \pi i (H_2 + iz^{-1})  }{k} \right) \right)}, &\mbox{if $2 \nmid k$},
 \end{array}
 \right.
\end{equation}
where $H_j$ is a solution to the congruence $jhH_j \equiv -1\pmod{k}$. 

Thus,
\begin{multline*}
\delta(n) = 
 i\underset{2\mid k}{\sum_{k=1}^{N}} k^{-1} \underset{(h,k)=1}
{ \sum_{0 \leqq h <  k}}  \frac{\omega(h,k)}{\omega(h,k/2)} e^{-2\pi i n h/k} 
\int_{z_I(h,k)}^{z_T(h,k)} 
\exp\left[  \frac{2\pi z }{k}
        \left(n+\frac{1}{24} \right)- \frac{\pi}{12kz}  \right] \\ \times
     F\left(  \exp \left( \frac{2\pi i (H_1 + iz^{-1})  }{k} \right) \right) \ dz \\
+ \frac{ i}{\sqrt{2}}\underset{2\nmid k}{\sum_{k=1}^{N} }k^{-1} \underset{(h,k)=1}
{ \sum_{0 \leqq h <  k}}  \frac{\omega(h,k)}{\omega(2h,k)} e^{-2\pi i n h/k} 
\int_{z_I(h,k)}^{z_T(h,k)} 
 \frac{\exp\left[  \frac{2\pi z }{k}
        \left(n+\frac{1}{24} \right) + \frac{\pi}{24 kz}  \right] }
    {F\left(  \exp \left( \frac{ \pi i (H_2 + iz^{-1})  }{k} \right) \right)} \ dz.
\end{multline*}
Next, we expand the appearances of $F$ as series:
\begin{multline*}
\delta(n) = 
 i\underset{2\mid k}{\sum_{k=1}^{N}} k^{-1} \underset{(h,k)=1}
{ \sum_{0 \leqq h <  k}}   \frac{\omega(h,k)}{\omega(h,k/2)} e^{-2\pi i n h/k} 
\int_{z_I(h,k)}^{z_T(h,k)} 
\exp\left[ 
        \frac{2\pi z }{k}
        \left(n+\frac{1}{24} \right)- \frac{\pi}{12kz}  \right] \\ \times
     \sum_{m=0}^\infty \delta(m) \exp \left( \frac{2\pi i (H_1 + iz^{-1}) m }{k} \right) \ dz \\
+ \frac{ i}{\sqrt{2}}\underset{2\nmid k}{\sum_{k=1}^{N} }k^{-1} \underset{(h,k)=1}
{ \sum_{0 \leqq h <  k}}   \frac{\omega(h,k)}{\omega(2h,k)}  e^{-2\pi i n h/k } 
\int_{z_I(h,k)}^{z_T(h,k)}
{\exp\left[   \frac{2\pi z }{k}
        \left(n+\frac{1}{24} \right) + \frac{\pi}{24 kz}  \right] }\\
   \times \sum_{m=0}^\infty \delta^{*}(m)
    { \exp \left( \frac{ \pi i (H_2 + iz^{-1})m  }{k} \right) } \ dz
\end{multline*} 
\begin{multline} \label{ThreeSum}
 = 
 i\underset{2\mid k}{\sum_{k=1}^{N}} 
 k^{-1}
  \sum_{m=0}^\infty \delta(m) \underset{(h,k)=1}
{ \sum_{0 \leqq h <  k}}   \frac{\omega(h,k)}{\omega(h,k/2)} 
 \exp \left[ \frac{2 \pi i}{k} (H_1m -hn  ) \right]
 \\ \times
\int_{z_I(h,k)}^{z_T(h,k)} 
\exp\left[ 
         \frac{2\pi z }{k}
        \left(n+\frac{1}{24} \right)- \frac{\pi}{kz} \left( 2m + \frac{1}{12} \right)  \right] \ dz \\
+ \frac{ i}{\sqrt{2}}\underset{2\nmid k}{\sum_{k=1}^{N} }k^{-1} 
 \sum_{m=1}^\infty \delta^{*}(m) \underset{(h,k)=1}
{ \sum_{0 \leqq h <  k}}   \frac{\omega(h,k)}{\omega(2h,k)} 
\exp \left[  \frac{\pi i}{k} (H_2 m - 2hn)  \right] \\ \times
\int_{z_I(h,k)}^{z_T(h,k)}
{\exp\left[  \frac{2\pi z }{k}
        \left(n+\frac{1}{24} \right) + \frac{\pi}{kz} \left( \frac{1}{24}- m \right)  \right] } \ dz \\
 + \frac{ i}{\sqrt{2}}\underset{2\nmid k}{\sum_{k=1}^{N} }k^{-1} 
 \underset{(h,k)=1}
{ \sum_{0 \leqq h <  k}}   \frac{\omega(h,k)}{\omega(2h,k)} 
 e^{-2\pi i n h/k  } 
\int_{z_I(h,k)}^{z_T(h,k)}
{\exp\left[  \frac{2\pi z }{k}
        \left(n+\frac{1}{24} \right) + \frac{\pi}{24kz}   \right] } \ dz
\end{multline} where 
\[ \frac{1}{F(q)} = \sum_{n=0}^\infty \delta^{*}(n) q^n ,\] and we have used the fact that
$\delta^{*}(0) = 1$.

 If the three sums in~\eqref{ThreeSum} are
designated $S_1$, $S_2$, and $S_3$
respectively, it can be shown via Kloosterman sum estimation
that $S_1, S_2 \to 0$ as $N\to\infty$ and only $S_3$ 
contributes to the final formula for $\delta(n)$.  
\begin{multline}
\frac{ -i}{24\sqrt{2}}\underset{2\nmid k}{\sum_{k=1}^{N} } \frac 1k
 \underset{(h,k)=1}
{ \sum_{0 \leqq h <  k}}   \frac{\omega(h,k)}{\omega(2h,k)} 
 e^{-2\pi i n h/k  } 
\int_{z_I(h,k)}^{z_T(h,k)}
{\exp\left[  \frac{2\pi z }{k}
        \left(n+\frac{1}{24} \right) + \frac{\pi}{24kz}   \right] } \ dz
\end{multline} 

Change variables
$t = \frac{\pi}{12k z}$ to obtain
\begin{align}
\delta(n) &= \frac{-i \pi}{24\sqrt{2} }
\underset{2\nmid k}{\sum_{k=1}^\infty} \frac{1}{k^2} \notag \\ &\quad\times
\underset{(h,k)=1}{\sum_{0\leqq h < k}} e^{-2\pi i h n/k}
\frac{ \omega(h,k)  }{\omega(2h,k) } 
 \int_{\pi/12-\infty i}^{\pi/12+\infty i} t^{-2}
  \exp\left( t + \frac{\pi^2 (24n+1)}{288k^2 t} \right) dt \notag \\
 &=  \frac{\pi}{\sqrt{24n+1} }
\underset{2\nmid k}{\sum_{k=1}^\infty} \frac 1k
\underset{(h,k)=1}{\sum_{0\leqq h < k}} e^{-2\pi i h n/k}
\frac{ \omega(h,k)  }{\omega(2h,k) } I_1\left( \frac{\pi\sqrt{24n+1}}{6k\sqrt{2}}\right).
\end{align}

\subsection{Summary of calculations}
We now summarize the required steps to find a Rademacher type formula for
$a(n)$ where 
\[ \sum_{n=0}^\infty a(n) q^n = \prod_{j=1}^J \frac{ f(q^{b_j})}
{f(q^{c_j})  } .\]

\begin{itemize}
  \item Find $L:=\mathrm{lcm}( b_1, b_2, \dots, b_J, c_1, c_2, \dots, c_J)$.
  \item For each divisor $d$ of $L$, there corresponds a case 
     $\gcd(k,L) = d$.
       \begin{itemize}
            \item To each case there corresponds to a summand of the form
            $$\Omega_{h,k} C \Psi_k(z) F(z,h,k)$$ which results from 
            applying the modular transformation~\eqref{EtaFunctionalEq} to that case.
            $\Omega_{h,k}$ is a product of powers of the $\omega$ 24$k$th root
            of $1$, $C$ is the constant that results, $\Psi_k(z)$ is the exponential
            expression, and $F(z,h,k)$ is the product of powers of $f$.
            \item Only those cases for which the co\"efficient of $z^{-1}$ in
               $\log \Psi_1(z)$ is positive 
               will contribute to the final formula; others can
               be shown to approach $0$ via Kloosterman sum estimation.
            \item Map $z \mapsto \pi/(12 k t)$.
            \item Evaluate integral in terms of the $I_1$ Bessel function.
       \end{itemize}
\end{itemize}



\section{Slater's list}
In 1952, L. J. Slater published a list of 130 identities of Rogers-Ramanujan
type~\cite{ljs1952}.  Many of the infinite products can be realized as
products of powers of $\eta$-functions, and have straightforward 
combinatorial interpretations as generating functions of restricted
classes of partitions or overpartitons.

Let us recall some of the identities in Slater's list.
\begin{align}
\sum_{n=0}^\infty \frac{(-1)^n q^{n(2n+1)}}{(q^2;q^2)_n(-q;q^2)_{n+1}}
&= \prod_{m=1}^\infty (1+q^{2m})(1-q^{2m-1}) \tag{S. 5}\\
  \sum_{n=0}^\infty \frac{ q^{n(n+1)/2} (-q)_n}{(q)_n} &=  
  \prod_{m=1}^\infty \frac{1-q^{4m}}{1-q^m} \tag{S. 8} \\
  \sum_{n=0}^\infty \frac{q^{n(2n+1)}}{(q)_{2n+1}} &=
  \prod_{m=1}^\infty (1+q^m) \tag{S. 9 = S. 84}\\
\sum_{n=0}^\infty \frac{q^{n^2} (-1)_{2n}}{(q^2;q^2)_n (q^2;q^4)_n} &= \prod_{m=1}^\infty
    \frac{1+q^{2m-1}}{1-q^{2m-1}}\tag{S. 10} \\
 \sum_{n=0}^\infty \frac{q^{n(n+1)} (-q;q^2)_{n}}{(q)_{2n+1} }
 &= \prod_{m=1}^\infty \frac{1-q^{4m}}{1-q^m} \tag{S. 11 = S.51 = S.64} \\
 \sum_{n=0}^\infty \frac{q^n (-1)_{2n}}{(q^2;q^2)_n} &= 
 \prod_{m=1}^\infty \frac{(1-q^{6m-3})^2 (1-q^{6m}) (1+q^m)}{1-q^m}
 \tag{S. 24}  \\
 \sum_{n=0}^\infty \frac{q^{n^2} (-q)_{n}}{(q;q^2)_{n+1} (q)_n} &= 
 \prod_{m=1}^\infty \frac{(1-q^{6m-3})^2 (1-q^{6m}) (1+q^m)}{1-q^m}
 \tag{S. 26}  \\
 \sum_{n=0}^\infty \frac{ q^{2n(n+1)} (-q;q^2)_n}{(q)_{2n+1}(-q^2;q^2)_n }
 &= \prod_{m=1}^\infty \frac{(1+q^{6m-5})(1+q^{6m-1})}
 {(1-q^{6m-4})(1-q^{6m-2})} \tag{S. 27}\\
   \sum_{n=0}^\infty \frac{q^{n(2n-1)}}{(q)_{2n}} 
  &=   \prod_{m=1}^\infty (1+q^m) \tag{S. 52 = S. 85}\\
   \sum_{n=0}^\infty \frac{q^{n(n+3)/2} (-q)_{n+1} (q^3;q^3)_n}
    {(q)_n (q)_{2n+2} } &= \prod_{m=1}^\infty
    \frac{(1-q^{18m})(1-q^{18m-3})(1-q^{18m-15})}{(1-q^{2m-1})(1-q^m)} 
    \tag{S. 76}\\
 \sum_{n=0}^\infty \frac{q^{n(n+1)/2} (-q)_n (q^3;q^3)_n}
    {(q)_n (q)_{2n+1} } &= \prod_{m=1}^\infty
    \frac{(1-q^{6m})(1+q^m) }{1-q^m} \tag{S. 77}\\
    1+\sum_{n=1}^\infty \frac{q^{n(n+1)/2} (-1)_{n+1} (q^3;q^3)_{n-1}}
    {(q)_{n-1} (q)_{2n} } &= \prod_{m=1}^\infty
    \frac{(1-q^{18m})(1-q^{18m-9})^2(1+q^m) }{1-q^m} \tag{S. 78}\\
    \sum_{n=0}^\infty \frac{q^{n(n+1)} (q^3;q^3)_n}{(q)_{2n+1} (q)_n}
    &= \prod_{m=1}^\infty \frac{1-q^{9m}}{1-q^m} \tag{S. 92} \\
    \sum_{n=0}^\infty \frac{ q^{n(n+1)} (q^3;q^6)_n (-q^2;q^2)_n}
    {(q^2;q^2)_{2n+1} (q;q^2)_n} 
   &=   \prod_{m=1}^\infty \frac{(1-q^{6m})(1+q^{12m-3})(1+q^{12m-9})
  }{ (1-q^{4m-2})(1-q^{2m})} \tag{S. 107} \\
    \sum_{n=0}^\infty \frac{ q^{n(n+2)} (q^3;q^6)_n (-q;q^2)_{n+1}}
    {(q^2;q^2)_{2n+1} (q;q^2)_n} &= \prod_{m=1}^\infty
    \frac{(1-q^{12m})(1-q^{4m-2})}{1-q^m}  \tag{S. 110 corrected}\\
    \sum_{n=0}^\infty \frac{q^{n(n+2)} (q^6;q^6)_n (-q;q^2)_{n+1}}
    {(q^2;q^2)_n  (q^2;q^2)_{2n+2} }
    &= \prod_{m=1}^\infty \frac{(1-q^{36m})(1-q^{36m-27})(1-q^{36m-9})}
    {(1-q^{2m-1})(1-q^{4m})}  \tag{S. 115}
\end{align}

Denote the c\"oefficient of $q^n$ in the power
series expansion of equation (S.$j$) above by
$S_j(n)$.
The following combinatorial interpretations are then immediate:

 \begin{itemize}
 \item $S_8(n)= S_{11}(n) = \delta_4(n)=$ the number of $4$-regular partitions of $n$; see Theorem~\ref{Hagis2}.
 \item $S_9(n) = S_{52}(n) = \delta(n)=$ the number of partitions into odd parts; 
  see Theorem~\ref{Hagis1}.
  \item $S_{10}(n) = $ the number of overpartitions of $n$ with only odd
  parts.
\item $S_{27}(n) = $ the number of overpartitions of $n$ where overlined
parts are odd nonmultiples of $3$ and the nonoverlined parts are even
nonmultiples of $6$.
 \item $S_{76}(n) =$ the number of overpartitions of $n$ where no
 nonoverlined part is congruent to 0, 3, or 15$\pmod{18}$.
\item $S_{77}(n)=$ the number of overpartitions of $n$ where no nonoverlined part is a multiple of $6$.   
\item $S_{92}(n)= \delta_9(n) = $ the number of $9$-regular partitions of $n$; 
see Theorem~\ref{Hagis2}.
\item $S_{107}(n) = $ the number of overpartitions of $n$ where
overlined parts are even or $\pm 3 \pmod{12}$ and nonoverlined
parts are $\pm 2\pmod{6}$.
\item $S_{110}(n) = $ the number of partitions of $n$ into parts not
congruent to $0,2,6,10 \pmod{12}$.
\item $S_{115}(n) = $ the number of partitions of $n$ into parts not
congruent to $0, \pm 9\pmod{36}$ nor congruent to $2\pmod{4}$.
\end{itemize}

The following Rademacher type formulas were conjectured with the
aid of \emph{Mathematica}:

\begin{multline} \label{s5}
S_5(n) = \frac{2\pi}{\sqrt{24n+1}}
 \underset{k\equiv 2\hskip -3mm\pmod{4}}{\sum_{k\geqq 1}} \frac 1k\\ \times
 \underset{(h,k)=1}{\sum_{0\leqq h<k}} e^{-2\pi i n h/k} 
 \frac{\omega(h,\frac k2)^2}{\omega(h,k) \omega(2h,\frac k2)}
 I_1 \left(  \frac{\pi\sqrt{24n+1}}{3k\sqrt{2}}   \right).
\end{multline}

 \begin{equation}
  S_{10}(n) = \frac{\pi}{4\sqrt{n}} \underset{2\nmid k}{\sum_{k\geqq 1}}
 \frac 1k \underset{(h,k)=1}{\sum_{0\leqq h<k}} e^{-2\pi i n h/k} 
  \frac{ \omega(h,k)^2 \omega(4h,k)}{\omega(2h,k)^3} 
  I_1\left(\frac{\pi\sqrt{n}}{k\sqrt{2}} \right).  \label{s10}
\end{equation}

\begin{equation} \label{s24}
  S_{24}(n) = \frac{\pi }{3  \sqrt{2n} }
  \underset{2,3\nmid k}{\sum_{k\geqq 1}}
  \frac{1}{k}   \underset{(h,k)=1}{\sum_{0\leqq h<k}} e^{-2\pi i n h/k} 
  \frac{   \omega(h,k)^2  \omega(6h,k )}
  {  \omega(2h,k) \omega(3h,k)^2 }
  I_1 \left(  \frac{\pi\sqrt{2n}}{k\sqrt{3}}  \right)
  \end{equation}
  

  \begin{multline} \label{s27}
S_{27}(n) = \frac{\pi }{ 9\sqrt{4n+1}} \sum_{d\mid 4} (d-2)(2d-5)
  \underset{(k,12)=d}{\sum_{k\geqq 1}}  \frac{1}{k }  \\ \times
   \underset{(h,k)=1}{\sum_{0\leqq h<k}} e^{-2\pi i n h/k}
   \frac{   \omega(h,k)  \omega( \frac{4h}{d}, \frac{k}{d} )  
    \omega( \frac{6h}{ \sqrt{d}}   ,  \frac{k}{\sqrt{d}} ) }
    { \omega( \frac{12 h}{d}, k  )    \omega( 3h  , \frac{k}{d}  ) 
      \omega( \frac{2h}{\sqrt{d}},  \frac{k}{\sqrt{d}}  ) } 
 I_1 \left(  \frac{\pi \sqrt{d(4n+1)} }{2k\sqrt{3}}  \right)
\end{multline}

\begin{multline} \label{s76}
  S_{76}(n) = \frac{\pi}{9\sqrt{2n+2}} \underset{(k,18)=1}{\sum_{k\geqq 1}}
 \frac 1k  \\ \times
 \underset{(h,k)=1}{\sum_{0\leqq h<k}} e^{-2\pi i n h/k} 
  \frac{ \omega(h,k)^2 \omega(6h,k) \omega(9h,k)}
  {\omega(2h,k) \omega(3h,k)  \omega(18h,k)^2} 
  I_1\left(\frac{2\pi\sqrt{2n+2}}{3k} \right)
\end{multline}

  \begin{equation}  \label{s77}
  S_{77}(n) = \frac{\pi\sqrt{2}}{3\sqrt{12n+3}} \underset{2,3\nmid k}
  {\sum_{k\geqq 1}}
 \frac 1k \underset{(h,k)=1}{\sum_{0\leqq h<k}} e^{-2\pi i n h/k} 
 \frac{ \omega(h,k)^2 }{\omega(2h,k) \omega(6h,k)} 
  I_1\left(\frac{\pi\sqrt{8n+2}}{3k} \right)
\end{equation}

 \begin{equation} \label{s78}
  S_{78}(n) = \frac{\pi\sqrt{2}}{9\sqrt{n}} \underset{(k,18)=1}{\sum_{k\geqq 1}}
 \frac 1k \underset{(h,k)=1}{\sum_{0\leqq h<k}} e^{-2\pi i n h/k} 
  \frac{ \omega(h,k)^2 \omega(18h,k) }{\omega(2h,k) \omega(9h,k)^2} 
  I_1\left(\frac{2\pi\sqrt{2n}}{3k} \right)
\end{equation}

\begin{multline} \label{s107}
S_{107}(n) = \frac{2\pi}{3\sqrt{24n+3}} \sum_{j=1}^2 \sqrt{4j-3}
\underset{(k,12)=j}{\sum_{k\geqq 1}} \frac 1k \\ \times
  \underset{(h,k)=1}{\sum_{0\leqq h<k}} e^{-2\pi i n h/k} 
 \frac{\omega(\frac{2h}{j}, \frac kj)^2  \omega(3h,k)
  \omega(\frac{12h}{j}, \frac{k}{j}) }
  {\omega(\frac{4h}{j}, \frac kj)  \omega(\frac{6h}{j}, \frac kj)^3}
  I_1 \left( \frac{\pi\sqrt{(3j-1)(8n+1)}}{6k} \right)
\end{multline}


\begin{multline}  \label{s110}
S_{110}(n) = \frac{2\pi}{9\sqrt{16n+6}}
\sum_{d\mid 4} \sqrt{(d-2)(7d-13)}
\underset{(k,12)=d}{\sum_{k\geqq 1}} 
 \frac 1k   \\ \times
 \underset{(h,k)=1}{\sum_{0\leqq h<k}} e^{-2\pi i n h/k} 
  \frac{ \omega(h,k) \omega( \frac{4h}{d}, \frac{k}{d}) }
  {\omega(\frac{2h}{\sqrt{d}}, \frac{k}{\sqrt{d}})
   \omega(\frac{12h}{d}, \frac{k}{d} ) }
  I_1\left(\frac{\pi\sqrt{1+d}\sqrt{8n+3}}{6k} \right)
\end{multline}

\begin{multline}  \label{s115}
S_{115}(n) = \frac{\pi}{27\sqrt{n+1}}
 \sum_{d\mid 4} (d-2)(2d-5)
 \underset{(k,36)=d}
{\sum_{k\geqq 1}}
\frac{1}{k}  
\\ \times\underset{(h,k)=1}{\sum_{0\leqq h<k}} e^{-2\pi i n h/k} 
\frac{ \omega(h,k) \omega(\frac{4h}{d},  \frac{k}{d} )
  \omega( \frac{18h}{\sqrt{d}}, \frac{k}{\sqrt{d}}  )  }
  { \omega(9h,k)  \omega( \frac{2h}{\sqrt{d}}  , \frac{k}{\sqrt{d}} )
  \omega(\frac{36h}{d},  \frac{k}{d} ) }
  I_1 \left(  \frac{2\sqrt{d}\pi \sqrt{n+1}}{3k} \right)
\end{multline}

\section{Numerical Test}
Each of the formulas \eqref{s10}--\eqref{s115}, along 
with Hagis's formula~\eqref{HagisDist} and Niven's formula
~\eqref{Niven}, was tested summing
$k$ from $1$ to $10$, and the value provided by the formula was
compared with the actual value.  In the chart below, the true value
of the given function at $n=100$ is provided along with the magnitude of the largest
error in the formula (when truncated at $k=10$) for $1\leqq n \leqq 100$.

\begin{center}
\begin{tabular}{| c | c | r |  r |}
\hline
Eq. no. &  function  &  value at $n=100$ & max error  \\
\hline
\eqref{HagisDist} & $\delta(n) $ & 444 793 & 0.211 \\
\eqref{Niven} & $S(n) $ & 20 901 & 0.318 \\
\eqref{s5} & $S_5(n)$  & 444 793 & 0.186 \\
\eqref{s10} & $S_{10}(n) $ & 29 025 326   & 0.210\\ 
\eqref{s24} & $S_{24}(n) $ & 793 378 722 & 0.200 \\
\eqref{s27} & $S_{27}(n) $ & 369 566  & 0.188 \\
\eqref{s76} & $S_{76}(n) $ & 15 008 235 468 & 0.050 \\
\eqref{s77} & $S_{77}(n) $ & 23 399 621 246 & 0.133 \\
\eqref{s78} & $S_{78}(n) $ & 26 086 456 322 & 0.143 \\
\eqref{s107} & $S_{107}(n) $ & 4 690 080 & 0.166 \\ 
\eqref{s110} & $S_{110}(n)$ & 4 731 983 & 0.216  \\
\eqref{s115} & $S_{115}(n)$  & 4 105 275 & 0.162 \\ \hline
\end{tabular}
\end{center}

\section*{Acknowledgments}
The author is grateful to the anonymous referee for carefully reading the manuscript, catching some typographical errors, and making a number of helpful suggestions.
Additionally, the author thanks Michael Schlosser for pointing out two errors that
had appeared in the 2010 published version in an email dated January 20, 2020.  The corrections are incorporated in this arXiv version of January 26, 2020.

\bibliographystyle{amsplain}

\end{document}